# REPRESENTATION OF FUNCTIONS IN SERIES WITH PARAMETER

## KHRISTO N. BOYADZHIEV

Department of mathematics, Ohio Northern University
Ada, Ohio, 45810, USA

E-mail: k-boyadzhiev@onu.edu

**Summary.** We prove a short general theorem which immediately implies some classical results of Hasse, Guillera and Sondow, Paolo Amore, and also Alzer and Richards. At the end we obtain a new representation for Euler's constant $\gamma$. The theorem transforms every Taylor series into a series depending on a parameter.

## 1 INTRODUCTION

In 1930 Hasse [1] proved a very interesting representation of Riemann's zeta function

$$\zeta(s) = \frac{1}{1-2^{1-s}} \sum_{n=0}^{\infty} \frac{1}{2^{n+1}} \sum_{k=0}^{n} \binom{n}{k} \frac{(-1)^k}{(k+1)^s} \qquad (1)$$

($\operatorname{Re} s > 1$). This results was rediscovered by Sondow [2] and later extended to the Lerch transcendent

$$\Phi(x,a,s) = \sum_{n=0}^{\infty} \frac{x^n}{(n+a)^s} \quad (\operatorname{Re} s > 1, a > 0, |x| \leq 1)$$

by Guillera and Sondow [3]. Namely, Guillera and Sondow showed that

$$\Phi(-x,a,s) = \frac{1}{x+1} \sum_{n=0}^{\infty} \left(\frac{x}{x+1}\right)^n \sum_{k=0}^{n} \binom{n}{k} \frac{(-1)^k}{(k+a)^s} \qquad (2)$$

for $\operatorname{Re} x > -1/2$. With $x = 1$ we have

$$\Phi(-1,1,s) = \sum_{n=0}^{\infty} \frac{(-1)^n}{(n+1)^s} = (1-2^{1-s})\zeta(s)$$

which shows that (2) extends (1).

Interesting representations of this nature were obtained by Amore [4], Alzer and Koumandos [5], and Alzer and Kendall [6]. For example, Amore [4] obtained the representation (extending a result of Flajolet and Vardi [7])

$$\pi = \sum_{n=1}^{\infty} \frac{\mu^n}{(\mu+1)^{n+1}} \left\{ \sum_{k=1}^{n} \binom{n}{k} \frac{1}{\mu^k} \frac{3^k-1}{4^k} \zeta(k+1) \right\} \qquad (3)$$

in a direct manner. He also showed that the parameter helps to improve the approximation of $\pi$ by the partial sums of the series in (3).

In this representations the series involves a free parameter, $\mu$. For example, among other things, Alzer and Kendall [6] obtained the interesting representation of the polylogarithm

$$\mathrm{Li}_s(x) = \frac{1}{\mu+1} \sum_{n=0}^{\infty} \left(\frac{\mu}{\mu+1}\right)^n \sum_{k=0}^{n} \binom{n}{k} \frac{x^{k+1}}{\mu^k (k+1)^s} \qquad (4)$$

where $s > 1, |x| \leq 1, \mu > -1/2$. This representation has the nice feature that the parameter does not appear on the left hand side and for a fixed $x$ one can change $\mu$. With $\mu = 1$ and $x \to -1$ the representation (4) implies (1) and with $x \to 1$ and it gives directly the representation

$$\zeta(s) = \sum_{n=0}^{\infty} \frac{1}{2^{n+1}} \sum_{k=0}^{n} \binom{n}{k} \frac{1}{(k+1)^s}.$$

The results in [6] are derived from an elaborate integral theorem. We present a more general theorem (with a simple proof) which implies immediately these results and many others. The theorem is true for any function $f(x)$ analytic in a neighborhood of the origin.

## 2 MAIN RESULTS

Given a function analytic in a neighborhood of the origin with the representation

$$f(x) = \sum_{n=0}^{\infty} a_n x^n$$

we have the following theorem.

**Theorem 1.** Suppose the function $f(t)$ defined by the series (1) is analytic in a disk $D = \{t, |t| < R\}$ containing the point $x$. Then we have the representation

$$f(\mu x) = \sum_{n=0}^{\infty} \frac{\mu^n}{(\mu+1)^{n+1}} \left\{ \sum_{k=0}^{n} \binom{n}{k} x^k a_k \right\} \qquad (5)$$

Where $\mu \in \mathbb{R}$ is an appropriate parameter with $-1/3 \leq \mu \leq 1$. Also,

$$f(x) = \sum_{n=0}^{\infty} \frac{1}{(\mu+1)^{n+1}} \left\{ \sum_{k=0}^{n} \binom{n}{k} \mu^{n-k} x^k a_k \right\}. \qquad (6)$$

As we will see later, the restriction $-1/3 \leq \mu \leq 1$ can be relaxed.

**Corollary 2.** With $\mu = 1$ in (5) we have

$$f(x) = \sum_{n=0}^{\infty} \frac{1}{2^{n+1}} \left\{ \sum_{k=0}^{n} \binom{n}{k} x^k a_k \right\}.$$

## 3. PROOF OF THE THEOREM

Let $L$ be a circle with radius $R_1 > |x|$, $R_1 < R$ such that $f(t) = a_0 + a_1 + ...$ is analytic in the disk $|t| < R_1$. For the coefficients $a_k$ we have

$$a_k = \frac{1}{2\pi i} \oint_L \frac{1}{\lambda^k} \frac{f(\lambda)}{\lambda} d\lambda,$$

Multiplying both sides by $\binom{n}{k} x^k$ and summing for $k = 0, 1, ..., n$ we find

$$\sum_{k=0}^{n} \binom{n}{k} x^k a_k = \frac{1}{2\pi i} \oint_L \left\{ \sum_{k=0}^{n} \binom{n}{k} \left(\frac{x}{\lambda}\right)^k \right\} \frac{f(\lambda)}{\lambda} d\lambda = \frac{1}{2\pi i} \oint_L \left(1 + \frac{x}{\lambda}\right)^n \frac{f(\lambda)}{\lambda} d\lambda .$$

Multiplying this by $\left(\frac{\mu}{\mu+1}\right)^n$ and summing for $n = 0, 1, ...$ we obtain

$$\sum_{n=0}^{\infty} \left(\frac{\mu}{\mu+1}\right)^n \sum_{k=0}^{n} \binom{n}{k} x^k a_k = \frac{1}{2\pi i} \oint_L \sum_{n=0}^{\infty} \left[\left(\frac{\mu}{\mu+1}\right)\left(1 + \frac{x}{\lambda}\right)\right]^n \frac{f(\lambda)}{\lambda} d\lambda .$$

The geometric series inside the integral has ratio $\frac{\lambda\mu + \mu x}{\lambda\mu + \lambda}$ and equals

$$\frac{1}{1 - \frac{\lambda\mu + \mu x}{\lambda\mu + \lambda}} = \frac{\lambda(\mu+1)}{\lambda - \mu x} \text{ provided } \left|\frac{\lambda\mu + \mu x}{\lambda\mu + \lambda}\right| < 1.$$

We integrate term by term and take out the terms $\left(\frac{\mu}{\mu+1}\right)^n$ to find

$$\frac{1}{2\pi i} \oint_L \sum_{n=0}^{\infty} \left[\left(\frac{\mu}{\mu+1}\right)\left(1 + \frac{x}{\lambda}\right)\right]^n \frac{f(\lambda)}{\lambda} d\lambda = \frac{\mu+1}{2\pi i} \oint_L \frac{f(\lambda)}{\lambda - \mu x} d\lambda = \sum_{n=0}^{\infty} \left(\frac{\mu}{\mu+1}\right)^n \sum_{k=0}^{n} \binom{n}{k} x^k a_k$$

$$= (\mu+1) f(\mu x).$$

From here

$$f(\mu x) = \sum_{n=0}^{\infty} \frac{\mu^n}{(\mu+1)^{n+1}} \sum_{k=0}^{n} \binom{n}{k} x^k a_k$$

and (5) is proved. For the convergence of the geometric series above we need

$$\left|\frac{\lambda\mu + \mu x}{\lambda\mu + \lambda}\right| < 1$$

and it is easy to see that this is true when $|\lambda| = R_1 > |x|$ and $-\frac{1}{3} \leq \mu \leq 1$. To show this we write

$$|\lambda\mu + \mu x| \leq |\lambda\mu| + |\mu x| = R_1|\mu| + |\mu||x| < 2R_1|\mu|$$

and we need $2R_1|\mu| \leq |\lambda|(\mu+1) = R_1(\mu+1)$, that is, $2|\mu| \leq \mu+1$, since we have already one strict inequality. If $\mu \geq 0$ this is $\mu \leq 1$. If $\mu < 0$, this is $-2\mu \leq \mu+1$, or $\mu \geq -1/3$.

To prove (6) we reason this way: as $|x| < R_1$ there is a region for $\mu$ where $\left|\frac{x}{\mu}\right| < R_1$. Then we replace $x$ by $x/\mu$ in (5) to get (6). By analytic continuation for $\mu$ we have (6).
The theorem is proved.

## 4 APPLICATIONS

For the beginning we will show how the theorem works with some elementary functions. Let first

$$f(x) = \frac{1}{1-x} = \sum_{n=0}^{\infty} x^n \quad (|x| < 1, \ a_n = 1, \ n = 0, 1, \ldots).$$

Then we have the representation

$$\frac{1}{1-x} = \sum_{n=0}^{\infty} \frac{\mu^n}{(\mu+1)^{n+1}} \left\{ \sum_{k=0}^{n} \binom{n}{k} \frac{x^k}{\mu^k} \right\}$$

which can be easily verified directly by using the formula for geometric series and the fact that

$$\sum_{k=0}^{n} \binom{n}{k} \frac{x^k}{\mu^k} = \left(1 + \frac{x}{\mu}\right)^n.$$

Replacing $x$ by $-x$ and integrating we find

$$\ln(x+1) = \sum_{n=0}^{\infty} \frac{\mu^n}{(\mu+1)^{n+1}} \left\{ \sum_{k=1}^{n} \binom{n}{k} \frac{(-1)^{k-1} x^k}{\mu^k k} \right\} \quad (|x| < 1).$$

If we take now $e^{-x} = \sum_{n=0}^{\infty} \frac{(-1)^n x^n}{n!}$ the theorem gives

$$e^{-x} = \sum_{n=0}^{\infty} \frac{\mu^n}{(\mu+1)^{n+1}} \left\{ \sum_{k=0}^{n} \binom{n}{k} \frac{(-x)^k}{\mu^k k!} \right\} = \sum_{n=0}^{\infty} \frac{\mu^n}{(\mu+1)^{n+1}} L_n\left(\frac{x}{\mu}\right) \qquad (7)$$

where $L_n(x)$ are the Laguerre polynomials. In particular, with $\mu = 1$

$$e^{-x} = \sum_{n=0}^{\infty} \frac{1}{2^{n+1}} L_n(x).$$

Taking into account the generating function for the Laguerre polynomials

$$\frac{1}{1-t} e^{-\frac{xt}{1-t}} = \sum_{n=0}^{\infty} L_n(x) t^n$$

the representation (7) can be verified directly.

For the binomial series we have the representation

$$(1+x)^{\beta} = \sum_{n=0}^{\infty} \frac{\mu^n}{(\mu+1)^{n+1}} \left\{ \sum_{k=0}^{n} \binom{n}{k} \binom{\beta}{k} \frac{x^k}{\mu^k} \right\} \quad (|x|<1).$$

Next we show that the result of Guillera and Sondow (2) as well as the result (4) of Alzer and Richards follow directly from the theorem.

**Corollary 3.** Equations (2) and (4) follow immediately from Theorem 1.

**Proof.**
The Lerch transcendent can be written as

$$\Phi(x,a,s) = \sum_{n=0}^{\infty} \frac{1}{(n+a)^s} x^n, \quad a_n = \frac{1}{(n+a)^s}.$$

Therefore, with $x=1$ and $\mu = x$ in (5) we obtain (2).

Note that in this case we need only $\mu > -\frac{1}{2}$ for the validity of the representation, as for such $\mu$ we have $\frac{|\mu|}{\mu+1} < 1$ and the series $\sum_{n=0}^{\infty} \left( \frac{\mu}{\mu+1} \right)^n$ converges which assures convergence of the series in (2) (see [3, p. 5]). The same remark is true for the representation (4) which we prove next. For the polylogarithm $\text{Li}_s(x) = \sum_{k=1}^{\infty} \frac{x^k}{k^s}$ we have

$$\frac{\text{Li}_s(x)}{x} = \sum_{n=0}^{\infty} \frac{x^n}{(n+1)^s}, \quad a_n = \frac{1}{(n+1)^s}$$

and the theorem gives

$$\frac{\text{Li}_s(x)}{x} = \sum_{n=0}^{\infty} \frac{\mu^n}{(\mu+1)^{n+1}} \left\{ \sum_{k=0}^{n} \binom{n}{k} \frac{x^k}{\mu^k} \frac{1}{(n+1)^s} \right\}$$

which is (4).

And now we proceed to reprove (3).

**Corollary 4**. The representation (3) is true

**Proof.**

We will use the Taylor series for the digamma function $\psi(x+1) = \dfrac{d}{dx}\Gamma(x+1)$

$$\psi(1-x) + \gamma = -\sum_{n=1}^{\infty} \zeta(n+1) x^n \qquad (8)$$

where $\gamma = -\psi(1)$ is Euler's constant, $\zeta(x)$ is Riemann's zeta function, and $|x| < 1$ [8, p. 159]. From the theorem

$$\psi(1-x) + \gamma = -\sum_{n=0}^{\infty} \frac{\mu^n}{(\mu+1)^{n+1}} \left\{ \sum_{k=0}^{n} \binom{n}{k} \frac{x^k}{\mu^k} \zeta(k+1) \right\}.$$

Next we write

$$\psi\left(1 - \frac{1}{4}\right) + \gamma = -\sum_{n=0}^{\infty} \frac{\mu^n}{(\mu+1)^{n+1}} \left\{ \sum_{k=0}^{n} \binom{n}{k} \frac{1}{\mu^k 4^k} \zeta(k+1) \right\}$$

$$\psi\left(\frac{1}{4}\right) + \gamma = \psi\left(1 - \frac{3}{4}\right) + \gamma = -\sum_{n=0}^{\infty} \frac{\mu^n}{(\mu+1)^{n+1}} \left\{ \sum_{k=0}^{n} \binom{n}{k} \frac{3^k}{\mu^k 4^k} \zeta(k+1) \right\}.$$

After subtraction

$$\psi\left(1 - \frac{1}{4}\right) - \psi\left(\frac{1}{4}\right) = \pi \cot\left(\frac{\pi}{4}\right) = \pi$$

by using the property $\psi(1-x) - \psi(x) = \pi \cot(\pi x)$. This proves (3).

In order to demonstrate the vast scope of Theorem 1 we will apply it now to the complete elliptic integrals of the first and second kind (see [9])

$$K(x) = \int_0^{\frac{\pi}{2}} \frac{d\theta}{\sqrt{1 - x^2 \sin^2 \theta}} = \int_0^1 \frac{dt}{\sqrt{(1-t^2)(1-x^2 t^2)}}$$

$$E(x) = \int_0^{\frac{\pi}{2}} \sqrt{1 - x^2 \sin^2 \theta}\ d\theta.$$

**Corollary 5.** The following representations hold

$$K(x) = \frac{\pi}{2} \sum_{n=0}^{\infty} \frac{\mu^n}{(\mu+1)^{n+1}} \left\{ \sum_{k=0}^{n} \binom{n}{k} \binom{2k}{k}^2 \frac{x^{2k}}{4^{2k} \mu^k} \right\} \qquad (9)$$

$$E(x) = \frac{\pi}{2} \sum_{n=0}^{\infty} \frac{\mu^n}{(\mu+1)^{n+1}} \left\{ \sum_{k=0}^{n} \binom{n}{k} \binom{2k}{k}^2 \frac{x^{2k}}{4^{2k}(1-2k)\mu^k} \right\}.$$

In particular, with $\mu = 1$

$$K(x) = \frac{\pi}{2} \sum_{n=0}^{\infty} \frac{1}{2^{n+1}} \left\{ \sum_{k=0}^{n} \binom{n}{k}\binom{2k}{k}^2 \frac{x^{2k}}{4^{2k}} \right\}$$

$$E(x) = \frac{\pi}{2} \sum_{n=0}^{\infty} \frac{1}{2^{n+1}} \left\{ \sum_{k=0}^{n} \binom{n}{k}\binom{2k}{k}^2 \frac{x^{2k}}{4^{2k}(1-2k)} \right\},$$

**Proof.**
The proof follows at once from the Taylor series for these functions

$$K(x) = \frac{\pi}{2} \sum_{n=0}^{\infty} \binom{2n}{n}^2 \frac{x^{2n}}{4^{2n}}, \quad E(x) = \frac{\pi}{2} \sum_{n=0}^{\infty} \binom{2n}{n}^2 \frac{x^{2n}}{4^{2n}(1-2n)}.$$

## 5. A SPECIAL CONSTANT

In [10] the author introduced the constant

$$M = \int_0^1 \frac{\psi(x+1)+\gamma}{x} dx \approx 1.257746$$

and proved, among other things, the representations

$$M = \sum_{n=1}^{\infty} \frac{(-1)^{n-1}\zeta(n+1)}{n} = \sum_{n=1}^{\infty} \frac{1}{n}\ln\left(1+\frac{1}{n}\right) = \sum_{n=1}^{\infty} \frac{\ln(n+1)}{n(n+1)} = \sum_{n=1}^{\infty} H_n\left(\zeta(n+1)-1\right)$$

$$= \sum_{n=1}^{\infty} \frac{1}{n}\left(n - \zeta(2) - \zeta(3) - \ldots - \zeta(n)\right)$$

with $H_n$ being the harmonic numbers. We give here one more representation of this constant depending on a parameter.

**Corollary 6**. For the constant $M$ we also have

$$M = \int_0^1 \frac{\psi(1+x)+\gamma}{x} dx = \sum_{n=0}^{\infty} \frac{\mu^n}{(\mu+1)^{n+1}} \left\{ \sum_{k=1}^{n} \binom{n}{k} \frac{(-1)^{k-1}}{\mu^k k} \zeta(k+1) \right\}.$$

**Proof**.
Changing $x$ to $-x$ in (8) we write

$$\psi(1+x)+\gamma = \sum_{n=1}^{\infty} (-1)^{n-1}\zeta(n+1)x^n$$

Here $a_0 = 0, a_n = (-1)^{n-1}\zeta(n+1)$ $(n>0)$. The theorem implies

$$\psi(1+x)+\gamma=\sum_{n=0}^{\infty}\frac{\mu^n}{(\mu+1)^{n+1}}\left\{\sum_{k=1}^{n}\binom{n}{k}\frac{x^k}{\mu^k}(-1)^{k-1}\zeta(k+1)\right\}. \qquad (10)$$

Here we divide both sides by $x$ and integrate from $0$ to $1$ to get

$$\int_0^1 \frac{\psi(1+x)+\gamma}{x}dx = \sum_{n=0}^{\infty}\frac{\mu^n}{(\mu+1)^{n+1}}\left\{\sum_{k=1}^{n}\binom{n}{k}\frac{(-1)^{k-1}}{\mu^k k}\zeta(k+1)\right\}$$

and the proof is finished.

In equation (10) above we have the parameter representation of the digamma function. For completeness we list here a similar representation for the log-gamma function. We just need to integrate both sides of (10) from $0$ to $x$ to get

$$\ln\Gamma(x+1)+\gamma x=\sum_{n=0}^{\infty}\frac{\mu^n}{(\mu+1)^{n+1}}\left\{\sum_{k=1}^{n}\binom{n}{k}\frac{x^{k+1}}{\mu^k(k+1)}(-1)^{k-1}\zeta(k+1)\right\}. \qquad (11)$$

From this representation we derive a parameter series for Euler's constant.

**Corollary 7.** For the Euler constant $\gamma$ we have the representation

$$\gamma=\sum_{n=0}^{\infty}\frac{\mu^n}{(\mu+1)^{n+1}}\left\{\sum_{k=1}^{n}\binom{n}{k}\frac{1}{\mu^k(k+1)}(-1)^{k-1}\zeta(k+1)\right\} \qquad (12)$$

and in particular, with $\mu=1$

$$\gamma=\sum_{n=0}^{\infty}\frac{1}{2^{n+1}}\left\{\sum_{k=1}^{n}\binom{n}{k}\frac{(-1)^{k-1}\zeta(k+1)}{(k+1)}\right\}.$$

Also,

$$\gamma=\frac{1}{\mu}(\mu+1-\ln(\mu+1))+\sum_{n=0}^{\infty}\frac{\mu^n}{(\mu+1)^{n+1}}\left\{\sum_{k=1}^{n}\binom{n}{k}\frac{(-1)^{k-1}}{\mu^k(k+1)}(\zeta(k+1)-1)\right\}. \qquad (13)$$

**Proof.**

Set $x=1$ in (10). The representation (12) follows immediately, as $\ln\Gamma(2)=\ln 1=0$. To prove (13) we write $\zeta(k+1)=\zeta(k+1)-1+1$ and use the binomial formula [11, p. 118]

$$\sum_{k=1}^{n}\binom{n}{k}(-1)^{k-1}\frac{1}{k+1}=\frac{n}{n+1}$$

to compute

$$\gamma=\sum_{n=0}^{\infty}\frac{\mu^n}{(\mu+1)^{n+1}}\left\{\sum_{k=1}^{n}\binom{n}{k}\frac{(-1)^{k-1}}{\mu^k(k+1)}(\zeta(k+1)-1)\right\}+\frac{1}{\mu}\sum_{n=0}^{\infty}\left(\frac{\mu}{\mu+1}\right)^{n+1}\frac{n}{n+1}$$

and (13) follows.

Alzer and Koumandos [5] obtained the representation

$$\gamma = \sum_{n=0}^{\infty} \frac{\mu^n}{(\mu+1)^{n+1}} \left\{ \sum_{k=0}^{n} \binom{n}{k} \frac{1}{\mu^k} (-1)^k S(k) \right\} \qquad (14)$$

where $S(k) = \sum_{n=1}^{\infty} \frac{1}{2^n + k}$. To compare (14) and (13) one can use the estimate from [10]

$$\frac{1}{2^{n+1}} < \zeta(n+1) - 1 < \frac{1}{2^{n+1}}(1 + \frac{2}{n}) \quad (n \geq 1).$$

## 6. CONCLUSIONS

We presented a general theorem for analytic functions which immediately implies some classical results for the Lerch transcendent $\Phi(x,a,s)$, the Riemann zeta function $\zeta(s)$, and the polylogarithm $\text{Li}_s(x)$. We also proved a new representation of Euler's constant $\gamma$ and a new representation of the constant $M$ introduced by the autor in [10].